# Lie Group analysis of Plasma equilibrium equations.

**Alexei F. Cheviakov**

**Queen's University at Kingston, 2002.**

Research advisor: Professor O. I. Bogoyavlenskij.

*In this paper we perform Lie group analysis of systems of partial differential equations which describe different cases of classical plasma equilibria, and find groups of transformations admitted by those equations in several important cases.*

*We also give a recipe how to apply Lie formalism to look for Bäcklund transformations between systems of partial differential equations, using plasma equilibrium system as an example.*

## *Contents.*

## 1. Lie group formalism.

Suppose we are given a set of equations
$$\mathbf{E}(\mathbf{u}, \mathbf{x}) = 0, \tag{1.0}$$
$\mathbf{E}$ can include differential operators. Let $\mathbf{u} = \mathbf{F}(\mathbf{x})$ be a solution of these equations. ($\mathbf{u} = \{u^1..u^m\}, \mathbf{x} = \{x^1..x^n\}$.)

The Lie group of transformations of the equations (1.0) is looked for in the form [1]
$$x_1^i = f^i(\mathbf{x},\mathbf{u},a), \tag{1.1}$$
$$u_1^k = g^k(\mathbf{x},\mathbf{u},a),$$
where the group operation is addition in $a$:
$$x_2^i = f^i(\mathbf{f}(\mathbf{x}_1,\mathbf{u}_1,a), \mathbf{g}(\mathbf{x}_1,\mathbf{u}_1,a),b) = f^i(\mathbf{x},\mathbf{u},a+b),$$
$$u_2^k = g^k(\mathbf{f}(\mathbf{x}_1,\mathbf{u}_1,a), \mathbf{g}(\mathbf{x}_1,\mathbf{u}_1,a),b) = g^k(\mathbf{x},\mathbf{u},a+b).$$
These variables $\mathbf{u}^1$ and $\mathbf{x}^1$ are supposed to satisfy the same equations as $\mathbf{u}$ and $\mathbf{x}$, i.e.
$$\mathbf{u}^1 = \mathbf{F}^1(\mathbf{x}^1)$$
must be a solution for the same differential equations (1.0) written in terms of $\mathbf{u}^1$ and $\mathbf{x}^1$.

To find the exact form of transformations, we search for the coordinates of vector field tangent to the orbits of transformation:
$$\xi^i(\mathbf{x},\mathbf{u}) = \left.\frac{\partial f^i(\mathbf{x},\mathbf{u},a)}{\partial a}\right|_{a=0}, \tag{1.2}$$
$$\eta^k(\mathbf{x},\mathbf{u}) = \left.\frac{\partial g^k(\mathbf{x},\mathbf{u},a)}{\partial a}\right|_{a=0}.$$

After we find these quantities, the Lie theorem lets us reconstruct the transformations - they are solutions to the equations
$$\frac{\partial f^i(a)}{\partial a} = \xi^i(\mathbf{f},\mathbf{g}), \quad \frac{\partial g^k(a)}{\partial a} = \eta^k(\mathbf{f},\mathbf{g}) \tag{1.3}$$
with initial conditions $f^i(0) = x^i$, $g^k(0) = u^k$.

To find the tangent vector field coordinates (1.2), we solve the determining equations
$$\left. \underset{1}{X}\, \mathbf{E}(\mathbf{u},\mathbf{x}) \right|_{\mathbf{E}(\mathbf{u},\mathbf{x})=0} = 0, \tag{1.4}$$
where $\underset{1}{X}$ is the infinitesimal operator of the prolonged vector field:
$$\underset{1}{X} \equiv \sum_i \xi^i \frac{\partial}{\partial x^i} + \sum_k \eta^k \frac{\partial}{\partial u^k} + \sum_{i,k} \xi_i^k \frac{\partial}{\partial u_i^k} = X + \sum_{i,k} \xi_i^k \frac{\partial}{\partial u_i^k}, \tag{1.5}$$
and $X$ is the infinitesimal operator.

Here $\xi_i^k$ are coordinates of prolonged tangent vector field corresponding to the derivatives $u_i^k$:
$$\xi_i^k(\mathbf{x},\mathbf{u}) = \left.\frac{\partial}{\partial a}\left(\frac{\partial u_1^k}{\partial x_1^i}\right)\right|_{a=0}. \tag{1.6}$$

These quantities are not unknown; they are connected with $\xi^i(\mathbf{x},\mathbf{u})$, $\eta^k(\mathbf{x},\mathbf{u})$ by

$$\xi_i^k(\mathbf{x},\mathbf{u}) = D_i\eta^k - u_m^k D_i\xi^m \qquad (1.7)$$

(summation over *m* is implied). Here $D_i$ is the *operator of complete differentiation*, defined as

$$D_i = \frac{\partial}{\partial x^i} + u_i^k \frac{\partial}{\partial u^k} + u_{ij}^k \frac{\partial}{\partial u_j^k} + \ldots \qquad (1.8)$$

Instead of transformations (1.1), we can also search for transformations in a richer class – namely, in the form

$$x_1^i = f^i(\mathbf{x}, \mathbf{u}, \underset{1}{\mathbf{u}}, a), \qquad (1.9)$$

$$u_1^k = g^k(\mathbf{x}, \mathbf{u}, \underset{1}{\mathbf{u}}, a),$$

where $\underset{1}{\mathbf{u}} = \left\{ \dfrac{\partial u^k}{\partial x^i},\ k=1..m,\ i=1..n \right\}$ are all derivatives.

## 2. Application of Lie group formalism to Plasma Equilibrium Equations.

Lie formalism has been applied to study different equations related to plasma physics and plasma confinement, for example, in [2], where the most general time-dependent system of equations describing motion of finitely-conducting plasma, Kadomtsev-Pogutse equations and some others are studied. But the Plasma Equilibrium equations, which are also of great importance to plasma physics, have not been studied in this way. It is known that simpler equations usually have richer groups of symmetries, - this gave us motivation to perform Lie analysis of this system and its important reductions.

### 2.1. Group analysis of Plasma Equilibrium equations.

We study the system of Plasma Equilibrium equations

$$\mathrm{curl}(\mathbf{B}) \times \mathbf{B} = \mathrm{grad}\, P, \qquad (2.1)$$

$$\mathrm{div}\, \mathbf{B} = 0.$$

Here **B** is magnetic field, *P* – plasma pressure.

The goal is to find all one-parametric groups of transformations of type (1.1) of the above system.

From here on we use the following notation for projections of the magnetic field:

**B** = (*A*, *B*, *C*).

The determining equations (1.4) demand that the differential equation $\underset{1}{X}\mathbf{E} = 0$ is solved under restriction $\mathbf{E}(\mathbf{u}, \mathbf{x}) = 0$, i.e. the original equations (2.1) must be satisfied at the same time. For that purpose, we express derivatives $P_x$, $P_y$, $P_z$, $C_z$ from equations (2.1)

$$C_z = -A_x - B_y,$$
$$P_x = CA_z - CC_x - BB_x + BA_y,$$
$$P_y = AB_x - AA_y - CC_y + CB_z,$$
$$P_z = BC_y - BB_z - AA_z + AC_x,$$
(2.2)

and use them solving the four equations $\underset{1}{X}\mathbf{E}(\mathbf{u}, \mathbf{x}) = 0$ for $\xi^i(x,y,z,A,B,C,P)$ and $\eta^k(x,y,z,A,B,C,P)$ ($i = 1..3$, $k = 1..4$). Note that these functions don't depend on the derivatives $\partial u^k/\partial x^i$, which are in turn connected only through (2.2), so all derivatives $\partial u^k/\partial x^i$ except $C_z$, $P_x$, $P_y$, $P_z$ are generally independent functions.

Therefore the method of solution is direct: put all coefficients at derivatives $A_x$, $A_y$, $A_z$, $B_x$, $B_y$, $B_z$, $C_x$, $C_y$ and their different products equal zero.

Let's call $XE_1..XE_3$ the results of the action of the $\underset{1}{X}$ operator on the projections of the first equation of (2.1); $XE_4$ – the result of the action of the $\underset{1}{X}$ operator on the second equation of (2.1) (Gauss theorem div $\mathbf{B} = 0$).

Collecting coefficients near $(\partial A/\partial x)^2$, $(\partial B/\partial x)^2$, $(\partial C/\partial x)^2$ in equations $XE_1..XE_4$, we find, first, that

$$\frac{\partial \xi^1}{\partial P} = 0,$$
(2.3)

and hence, by symmetry of coordinates, $\dfrac{\partial \xi^2}{\partial P} = \dfrac{\partial \xi^3}{\partial P} = 0$.

We also find that the following takes place:

$$\frac{\partial \xi^1}{\partial C} = \frac{\partial \xi^3}{\partial A},\ B\frac{\partial \xi^2}{\partial P} = \frac{\partial \xi^2}{\partial B} = 0,\ C\frac{\partial \xi^3}{\partial P} = \frac{\partial \xi^3}{\partial C} = 0,$$
$$\frac{\partial \xi^1}{\partial B} = A\frac{\partial \xi^2}{\partial P} = 0,\ \frac{\partial \xi^1}{\partial C} = A\frac{\partial \xi^3}{\partial P} = 0,\ \frac{\partial \xi^3}{\partial C} = -\frac{\partial \xi^1}{\partial A} = 0.$$
(2.4)

Therefore, using (2.3) and its consequence, we conclude that

$$\frac{\partial \xi^1}{\partial A} = 0,\ \frac{\partial \xi^1}{\partial B} = 0,\ \frac{\partial \xi^1}{\partial C} = 0,\ \frac{\partial \xi^2}{\partial B} = 0,\ \frac{\partial \xi^3}{\partial A} = 0,\ \frac{\partial \xi^3}{\partial C} = 0.$$
(2.5)

Using (2.3), (2.5) and calculating coefficients of $(\partial A/\partial y)^2$, $(\partial B/\partial y)^2$, $(\partial C/\partial y)^2$ in equations $XE_1..XE_4$, we also observe that

$$\frac{\partial \xi^2}{\partial A} = 0,\ \frac{\partial \xi^2}{\partial C} = 0,\ \frac{\partial \xi^3}{\partial B} = 0,$$
(2.6)

which lets us conclude that the coordinates of tangent vector field corresponding to spatial variables $x, y, z$ depend only on them, and not on the plasma parameters **B**, $P$:

$$\xi^i(\mathbf{x}, \mathbf{u}) = \xi^i(\mathbf{x}). \tag{2.7}$$

Now we substitute (2.7) into the equations $XE_1..XE_3$ and set coefficients at derivatives $A_x$, $A_y$, $A_z$, $B_x$, $B_y$, $B_z$, $C_x$, $C_y$ to zero. This results in the following equalities:

$$\frac{\partial \eta^1}{\partial B} = -\frac{\partial \xi^2}{\partial x}, \quad \frac{\partial \eta^1}{\partial C} = -\frac{\partial \xi^3}{\partial x}, \quad \frac{\partial \eta^2}{\partial A} = -\frac{\partial \xi^1}{\partial y},$$

$$\frac{\partial \eta^2}{\partial C} = -\frac{\partial \xi^3}{\partial y}, \quad \frac{\partial \eta^3}{\partial B} = -\frac{\partial \xi^2}{\partial z}, \quad \frac{\partial \eta^3}{\partial A} = -\frac{\partial \xi^1}{\partial z}, \tag{2.8}$$

$$\frac{\partial \eta^4}{\partial A} = \frac{\partial \eta^4}{\partial B} = \frac{\partial \eta^4}{\partial C} = 0.$$

The last three equalities mean that the tangent vector field coordinate corresponding to pressure $\eta^4$ is independent of the components of the magnetic field $A, B, C$:

$$\eta^4 = \eta^4(x, y, z, P). \tag{2.9}$$

Adding and subtracting the equations of (2.8) yields the following relations:

$$\frac{\partial \xi^1}{\partial y} = -\frac{\partial \xi^2}{\partial x}, \quad \frac{\partial \xi^1}{\partial z} = -\frac{\partial \xi^3}{\partial x}, \quad \frac{\partial \xi^2}{\partial z} = -\frac{\partial \xi^3}{\partial y},$$

$$\frac{\partial \eta^1}{\partial P} B = \frac{\partial \eta^2}{\partial P} A, \quad \frac{\partial \eta^3}{\partial P} A = \frac{\partial \eta^1}{\partial P} C, \quad \frac{\partial \eta^2}{\partial P} C = \frac{\partial \eta^3}{\partial P} B,$$

i.e.

$$\frac{\partial \xi^i}{\partial x^j} = -\frac{\partial \xi^j}{\partial x^i}, i \neq j; \quad \frac{\partial \eta^k}{\partial P} u^m = \frac{\partial \eta^m}{\partial P} u^k, k \neq m. \tag{2.10}$$

Using (2.7) and (2.10), from (2.8) we conclude that tangent vector field components corresponding to projections of **B** must have the form:

$$\eta^1 = B \cdot \xi^1_y + C \cdot \xi^1_z + f^A(A, P, x, y, z)$$

$$\eta^2 = -A \cdot \xi^1_y + C \cdot \xi^2_z + f^B(B, P, x, y, z) \tag{2.11}$$

$$\eta^3 = -A \cdot \xi^1_z - B \cdot \xi^2_z + f^C(C, P, x, y, z)$$

Applying relations (2.10) to formulas (2.11), we obtain

$$B \frac{\partial f^A}{\partial P} = A \frac{\partial f^B}{\partial P}, \quad A \frac{\partial f^C}{\partial P} = C \frac{\partial f^A}{\partial P}, \quad C \frac{\partial f^B}{\partial P} = B \frac{\partial f^C}{\partial P},$$

which means

$$\frac{\partial}{\partial P} \frac{f^A}{A} = \frac{\partial}{\partial P} \frac{f^B}{B} = \frac{\partial}{\partial P} \frac{f^C}{C} = \frac{\partial \varphi}{\partial P}(x, y, z, P). \tag{2.12}$$

Taking into consideration the lists of arguments of $f^A, f^B, f^C$ (see (2.11)), we find that these functions have the form

$$f^A = A \cdot \varphi(x, y, z, P) + \omega^A(x, y, z, A),$$
$$f^B = B \cdot \varphi(x, y, z, P) + \omega^B(x, y, z, B), \quad (2.13)$$
$$f^C = C \cdot \varphi(x, y, z, P) + \omega^C(x, y, z, C).$$

Now return to the determining equation $XE_4$. The coefficients at the derivatives $A_x, B_y$ there are

$$\frac{\partial \xi^3}{\partial z} - \frac{\partial \eta^3}{\partial C} - \frac{\partial \xi^1}{\partial x} + \frac{\partial \eta^1}{\partial A} = 0, \quad \frac{\partial \xi^3}{\partial z} - \frac{\partial \eta^3}{\partial C} - \frac{\partial \xi^2}{\partial y} + \frac{\partial \eta^2}{\partial B} = 0. \quad (2.14)$$

After substituting here the form of $\eta^k$ that we have found in (2.11), we find that

$$\frac{\partial}{\partial A}\omega^A(x,y,z,A), \quad \frac{\partial}{\partial B}\omega^B(x,y,z,B), \quad \frac{\partial}{\partial C}\omega^C(x,y,z,C)$$

are functions only of $x, y, z$. Therefore

$$f^A(A, P, x, y, z) = A \cdot \left(\varphi(x, y, z, P) + c^A(x, y, z)\right) + w^A(x, y, z),$$
$$f^B(B, P, x, y, z) = B \cdot \left(\varphi(x, y, z, P) + c^B(x, y, z)\right) + w^B(x, y, z), \quad (2.15)$$
$$f^C(C, P, x, y, z) = C \cdot \left(\varphi(x, y, z, P) + c^C(x, y, z)\right) + w^C(x, y, z).$$

Now we substitute this into the free-of-derivatives part of the equation $XE_3$

$$B\frac{\partial \eta^3}{\partial y} - B\frac{\partial \eta^2}{\partial z} - \frac{\partial \eta^4}{\partial z} - A\frac{\partial \eta^1}{\partial z} + A\frac{\partial \eta^3}{\partial x} = 0$$

and collect coefficients near $A^2$, to find that $\frac{\partial}{\partial z}\varphi(x, y, z, P)$ is independent of $P$. Hence from general equivalence of $x, y, z,$ we get

$$\varphi(x, y, z, P) = \Phi(x, y, z) + \theta(P), \quad (2.16)$$

and (2.15) can be rewritten as

$$f^A(A, P, x, y, z) = A \cdot \left(q^A(x, y, z) + \theta(p)\right) + w^A(x, y, z),$$
$$f^B(B, P, x, y, z) = B \cdot \left(q^B(x, y, z) + \theta(p)\right) + w^C(x, y, z), \quad (2.17)$$
$$f^C(C, P, x, y, z) = C \cdot \left(q^C(x, y, z) + \theta(p)\right) + w^C(x, y, z).$$

Plugging the functions $\eta^k$ that are found from (2.11) and (2.17) into free-of-derivatives parts of $XE_1, XE_2, XE_3$ and collecting terms at $AC, AB$ and $BC$ respectively, we get

$$\frac{\partial^2 \xi^1}{\partial x \partial z} = -\frac{\partial q^A}{\partial z}, \quad \frac{\partial^2 \xi^1}{\partial y^2} = -\frac{\partial q^B}{\partial x}, \quad \frac{\partial^2 \xi^2}{\partial z^2} = -\frac{\partial q^C}{\partial y}. \quad (2.18)$$

Using these formulas and relations (2.10), we get

$$q^A = -\frac{\partial \xi^1}{\partial x} + g_1(x, y), \quad q^B = -\frac{\partial \xi^2}{\partial y} + g_2(y, z), \quad q^C = -\frac{\partial \xi^3}{\partial z} + g_3(x, z) \quad (2.19)$$

But collecting terms at $C^2$, $A^2$ and $B^2$ respectively, we get

$$\frac{\partial^2 \xi^1}{\partial z^2} = \frac{\partial q^C}{\partial x}, \quad \frac{\partial^2 \xi^1}{\partial x \partial y} = -\frac{\partial q^A}{\partial y}, \quad \frac{\partial^2 \xi^2}{\partial y \partial z} = -\frac{\partial q^B}{\partial z}, \quad (2.20)$$

so that

$$\frac{\partial}{\partial y} g_1(x, y) = 0, \quad \frac{\partial}{\partial z} g_2(y, z) = 0, \quad \frac{\partial}{\partial x} g_3(x, z) = 0.$$

Hence

$$q^A = -\frac{\partial \xi^1}{\partial x} + g_1(x), \quad q^B = -\frac{\partial \xi^2}{\partial y} + g_2(y), \quad q^C = -\frac{\partial \xi^3}{\partial z} + g_3(z) \quad (2.21)$$

If in free-of-derivatives parts of $XE_1$, $XE_2$, $XE_3$ we collect terms independent of $A$, $B$, $C$, we will get

$$\frac{\partial \eta^4}{\partial x} = \frac{\partial \eta^4}{\partial y} = \frac{\partial \eta^4}{\partial z} = 0,$$

so

$$\eta^4 = \eta^4(P). \quad (2.22)$$

Substituting this new result into $XE_1$ and collecting terms near $C \cdot A^3 \cdot C_x$, we find

$$\frac{\partial}{\partial P} \theta(P) = 0. \quad (2.23)$$

and (2.17) rewrites as

$$f^A(A, P, x, y, z) = A \cdot \left(-\xi_x^x + g_1(x)\right) + w^A(x, y, z),$$
$$f^B(B, P, x, y, z) = B \cdot \left(-\xi_y^y + g_2(y)\right) + w^B(x, y, z), \quad (2.24)$$
$$f^C(C, P, x, y, z) = C \cdot \left(-\xi_z^z + g_3(z)\right) + w^C(x, y, z).$$

Substituting (2.11), (2.24) into $XE_1$ and collecting terms near $C_x$, we obtain

$$w^C(x, y, z) = 0$$

and

$$2 \cdot q^C = 2 \cdot \left(-\frac{\partial \xi^3}{\partial z} + g_3(z)\right) = \frac{\partial \eta^4(P)}{\partial P}.$$

The first relation lets us conclude, by symmetry of notation, that

$$w^A(x, y, z) = w^B(x, y, z) = w^C(x, y, z) = 0, \quad (2.25)$$

and the second – that

$$\eta^4(P) = 2\alpha P + \beta, \quad (2.26)$$

and again by symmetry of coordinates

$$\alpha = q^A = q^B = q^C = \text{const.} \quad (2.27)$$

(here $\alpha = \text{const}$, $\beta = \text{const}$).

Altogether, the tangent vector field coordinates corresponding to **B** and $P$ are

$$\eta^1 = \alpha \cdot A + B \cdot \xi^1_y + C \cdot \xi^1_z,$$
$$\eta^2 = \alpha \cdot B - A \cdot \xi^1_y + C \cdot \xi^2_z, \quad (2.28)$$
$$\eta^3 = \alpha \cdot C - A \cdot \xi^1_z - B \cdot \xi^2_z,$$
$$\eta^4 = 2\alpha P + \beta.$$

Using relations (2.14), (2.18), (2.20), (2.27), (2.28) together, we find the following conditions relations on $\xi^i$:

$$\frac{\partial^2 \xi^1}{\partial y^2} = 0,\ \frac{\partial^2 \xi^1}{\partial z^2} = 0,\ \frac{\partial^2 \xi^1}{\partial x \partial y} = 0,\ \frac{\partial^2 \xi^1}{\partial x \partial z} = 0,\ \frac{\partial^2 \xi^2}{\partial z^2} = 0,\ \frac{\partial^2 \xi^2}{\partial y \partial z} = 0,$$

$$\frac{\partial^2 \xi^1}{\partial y \partial z} = 0,\ \frac{\partial^2 \xi^2}{\partial x \partial z} = 0,\ \frac{\partial \xi^1}{\partial x} = \frac{\partial \xi^2}{\partial y} = \frac{\partial \xi^3}{\partial z}, \quad (2.29)$$

After some calculations, it's easy to see that $\xi^1, \xi^2, \xi^3$ must be linear. Making use of (2.29), we write them as

$$\xi^1 = R \cdot x + a_2 \cdot y + a_3 \cdot z + a_4,$$
$$\xi^2 = b_1 \cdot x + R \cdot y + b_3 \cdot z + b_4, \quad (2.30)$$
$$\xi^3 = c_1 \cdot x + c_2 \cdot y + R \cdot z + c_4$$

From the first formula of (2.11), we get the following connections for coefficients in (2.30):

$$a_2 = -b_1 = \delta,\ a_3 = -c_1 = \gamma,\ b_3 = -c_2 = \mu.$$

Finally, the tangent vector field coordinates corresponding to spatial variables are

$$\xi^1 = R \cdot x + \delta \cdot y + \gamma \cdot z + l_1,$$
$$\xi^2 = -\delta \cdot x + R \cdot y + \mu \cdot z + l_2, \quad (2.31)$$
$$\xi^3 = -\gamma \cdot x - \mu \cdot y + R \cdot z + l_3,$$

The tangent vector field coordinates corresponding to plasma parameters are

$$\eta^1 = \alpha \cdot A + \delta \cdot B + \gamma \cdot C,$$
$$\eta^2 = -\delta \cdot A + \alpha \cdot B + \mu \cdot C, \quad (2.32)$$
$$\eta^3 = -\gamma \cdot A - \mu \cdot B + \alpha \cdot C,$$
$$\eta^4 = 2 \cdot \alpha \cdot P + \beta.$$

In (2.31), (2.32) $\alpha, \beta, \delta, \gamma, \mu, R, l_1, l_2, l_3$ are arbitrary constants.

Altogether, the complete expression for the infinitesimal operator $X$ admitted by (2.1) is

$$X = \left(l_1 \frac{\partial}{\partial x} + l_2 \frac{\partial}{\partial y} + l_3 \frac{\partial}{\partial z} + \beta \frac{\partial}{\partial P}\right)$$
$$+ R \cdot \left(x \frac{\partial}{\partial x} + y \frac{\partial}{\partial y} + z \frac{\partial}{\partial z}\right) + \alpha \cdot \left(A \frac{\partial}{\partial A} + B \frac{\partial}{\partial B} + C \frac{\partial}{\partial C} + 2P \frac{\partial}{\partial P}\right)$$
$$+ \delta \cdot \left(y \frac{\partial}{\partial x} - x \frac{\partial}{\partial y} + B \frac{\partial}{\partial A} - A \frac{\partial}{\partial B}\right) + \gamma \cdot \left(z \frac{\partial}{\partial x} - x \frac{\partial}{\partial z} + C \frac{\partial}{\partial A} - A \frac{\partial}{\partial C}\right)$$
$$+ \mu \cdot \left(z \frac{\partial}{\partial y} - y \frac{\partial}{\partial z} + C \frac{\partial}{\partial B} - B \frac{\partial}{\partial C}\right)$$

(2.33)

## 2.2. Discussion and summary.

Let us analyze the operator (2.33). All the bracketed terms in it are independent, because they depend on different arbitrary constants.

The first term corresponds to four independent one-parametric shifts $G_{sh} \times G_{sh} \times G_{sh} \times G_{sh}$:

$$x' = x + a \cdot l_1,\ y' = y + a \cdot l_2,\ z' = z + a \cdot l_3,\ P' = P + a \cdot \beta,$$ (2.34)

which can be easily verified directly using (1.2), (1.3).

The second term describes two independent scaling transformations $G_{sc} \times G_{sc}$:

$$x' = e^R \cdot x,\ y' = e^R \cdot y,\ z' = e^R \cdot z,$$
$$A' = e^\alpha \cdot A,\ B' = e^\alpha \cdot B,\ C' = e^\alpha \cdot C,$$
$$P' = e^{2\alpha} \cdot P.$$

(2.35)

The other three terms are rotations in planes XY, XZ and YZ respectively. For example, the third term of (2.33) (with $\delta$) corresponds to the rotation transformation

$$x' = x \cdot \cos(a\delta) + y \cdot \sin(a\delta),$$
$$y' = y \cdot \cos(a\delta) - x \cdot \sin(a\delta),$$
$$z' = z,$$
$$A' = A \cdot \cos(a\delta) + B \cdot \sin(a\delta),$$
$$B' = B \cdot \cos(a\delta) - A \cdot \sin(a\delta),$$
$$C' = C,\ P' = P.$$

(2.36)

Similar formulas express rotations in XZ and YZ planes. Together all three rotations form an SO(3) group.

The groups listed above make up the general group of transformations of type (1.1) admissible by the system (2.1):

$$G_{TOT} = (G_{sh} \times G_{sh} \times G_{sh} \times G_{sh}) \times (G_{sc} \times G_{sc}) \times SO(3).$$ (2.37)

<u>Conclusion</u>: the whole group of transformations admitted by (2.1) consists of space and pressure shifts, scaling of space coordinates, magnetic field and pressure, and space rotations SO(3).

## 3. Group analysis of Plasma Equilibrium Equations up to 1<sup>st</sup> derivatives.

Now we extend our study of the Plasma Equilibrium system (2.1) to the case where Lie point transformations depend not only on dependent and independent variables, but also on all derivatives, i.e. repeat the procedure described in Section 2 using transformations (1.9) instead of (1.1).

The coordinates of tangent vector fields will then have the form

$$\xi^i = \xi^i(\mathbf{x},\mathbf{u},\underset{1}{\mathbf{u}}), \quad \eta^k = \eta^k(\mathbf{x},\mathbf{u},\underset{1}{\mathbf{u}}).\tag{3.1}$$

First we find the coordinates of the prolonged tangent vector field corresponding to derivatives (1.7), using (1.8), and then solve equations (1.4), setting coefficients near all *second derivatives* to zero.

After involved computations, the result appears to be exactly the same as in Section 2, i.e.

$$\begin{aligned}X =\;& \left(l_1\frac{\partial}{\partial x}+l_2\frac{\partial}{\partial y}+l_3\frac{\partial}{\partial z}+\beta\frac{\partial}{\partial P}\right)\\&+R\cdot\left(x\frac{\partial}{\partial x}+y\frac{\partial}{\partial y}+z\frac{\partial}{\partial z}\right)+\alpha\cdot\left(A\frac{\partial}{\partial A}+B\frac{\partial}{\partial B}+C\frac{\partial}{\partial C}+2P\frac{\partial}{\partial P}\right)\\&+\delta\cdot\left(y\frac{\partial}{\partial x}-x\frac{\partial}{\partial y}+B\frac{\partial}{\partial A}-A\frac{\partial}{\partial B}\right)+\gamma\cdot\left(z\frac{\partial}{\partial x}-x\frac{\partial}{\partial z}+C\frac{\partial}{\partial A}-A\frac{\partial}{\partial C}\right)\\&+\mu\cdot\left(z\frac{\partial}{\partial y}-y\frac{\partial}{\partial z}+C\frac{\partial}{\partial B}-B\frac{\partial}{\partial C}\right)\end{aligned}\tag{3.2}$$

<u>Conclusion:</u> the Plasma Equilibrium Equations (2.1) do not have any Lie transformations that depend on derivatives.

## 4. Group analysis of pure magnetic field equations (up to 1<sup>st</sup> derivatives).

Let us study groups of transformations admitted by a magnetic field in vacuum, which is a special subclass of force-free magnetic fields.

$$\text{curl}\,\mathbf{B}=0,\tag{4.1}$$
$$\text{div}\,\mathbf{B}=0.$$

We looked for transformations that depend on the magnetic field and its derivatives, but not on the coordinates:

$$\xi^i = \xi^i(\mathbf{u},\underset{1}{\mathbf{u}}), \quad \eta^k = \eta^k(\mathbf{u},\underset{1}{\mathbf{u}}).$$

When the transformations depend also on the coordinates, the determining equations are much harder to solve.

The solution of the equations (1.4) then gives us the infinitesimal operator admitted by (4.1):

$$X = f_1(\mathbf{u},\underset{1}{\mathbf{u}})\frac{\partial}{\partial x} + f_2(\mathbf{u},\underset{1}{\mathbf{u}})\frac{\partial}{\partial y} + f_3(\mathbf{u},\underset{1}{\mathbf{u}})\frac{\partial}{\partial z}$$
$$+ \left(\frac{\partial A}{\partial x}\cdot f_1(\mathbf{u},\underset{1}{\mathbf{u}}) + \frac{\partial B}{\partial x}\cdot f_2(\mathbf{u},\underset{1}{\mathbf{u}}) + \frac{\partial C}{\partial x}\cdot f_3(\mathbf{u},\underset{1}{\mathbf{u}}) + \alpha\cdot A + \beta\cdot\frac{\partial A}{\partial x} + \gamma\cdot\frac{\partial B}{\partial x} + \delta\cdot\frac{\partial C}{\partial x} + \mu_1\right)\frac{\partial}{\partial A}$$
$$+ \left(\frac{\partial B}{\partial x}\cdot f_1(\mathbf{u},\underset{1}{\mathbf{u}}) + \frac{\partial B}{\partial y}\cdot f_2(\mathbf{u},\underset{1}{\mathbf{u}}) + \frac{\partial C}{\partial y}\cdot f_3(\mathbf{u},\underset{1}{\mathbf{u}}) + \alpha\cdot B + \beta\cdot\frac{\partial B}{\partial x} + \gamma\cdot\frac{\partial B}{\partial y} + \delta\cdot\frac{\partial C}{\partial y} + \mu_2\right)\frac{\partial}{\partial B} \quad (3.2)$$
$$+ \left(\frac{\partial C}{\partial x}\cdot f_1(\mathbf{u},\underset{1}{\mathbf{u}}) + \frac{\partial C}{\partial x}\cdot f_2(\mathbf{u},\underset{1}{\mathbf{u}}) + \frac{\partial C}{\partial z}\cdot f_3(\mathbf{u},\underset{1}{\mathbf{u}}) + \alpha\cdot C + \beta\cdot\frac{\partial C}{\partial x} + \gamma\cdot\frac{\partial C}{\partial y} + \delta\cdot\frac{\partial C}{\partial z} + \mu_3\right)\frac{\partial}{\partial C}.$$

This group of transformations is richer than (2.33) or (3.2) for PE equations, and more interesting as it does depend explicitly on first derivatives.

## 5. Group analysis of force-free plasma equilibrium equations (up to 1$^{st}$ derivatives).

In this chapter we look for Lie-Bäcklund transformations admissible by the system of force-free plasma equilibrium equations

$$\operatorname{curl}\mathbf{B} = \alpha(\mathbf{x})\mathbf{B}, \quad (5.1)$$
$$\operatorname{div}\mathbf{B} = 0,$$

(We again use notation $\mathbf{B} = (A, B, C)$.)

Let us allow the transformed variables to depend on not only on original variables, but also on derivatives (see Section 1):

$$x_1^i = f^i(\mathbf{x},\mathbf{u},\underset{1}{\mathbf{u}},a), \quad (5.2)$$
$$u_1^k = g^k(\mathbf{x},\mathbf{u},\underset{1}{\mathbf{u}},a).$$

Repeating the procedure performed in Section 2 for PE equations, after involved computations we find that the system (5.1) admits the infinitesimal operator

$$X = \left(l_1\frac{\partial}{\partial x} + l_2\frac{\partial}{\partial y} + l_3\frac{\partial}{\partial z}\right)$$
$$+ R\cdot\left(x\frac{\partial}{\partial x} + y\frac{\partial}{\partial y} + z\frac{\partial}{\partial z}\right) + \alpha\cdot\left(A\frac{\partial}{\partial A} + B\frac{\partial}{\partial B} + C\frac{\partial}{\partial C}\right)$$
$$+ \delta\cdot\left(y\frac{\partial}{\partial x} - x\frac{\partial}{\partial y} + B\frac{\partial}{\partial A} - A\frac{\partial}{\partial B}\right) + \gamma\cdot\left(z\frac{\partial}{\partial x} - x\frac{\partial}{\partial z} + C\frac{\partial}{\partial A} - A\frac{\partial}{\partial C}\right) \quad (5.3)$$
$$+ \mu\cdot\left(z\frac{\partial}{\partial y} - y\frac{\partial}{\partial z} + C\frac{\partial}{\partial B} - B\frac{\partial}{\partial C}\right)$$

which *does not* depend on derivatives. This operator coincides with the operator (3.2) for Plasma Equilibrium equations (with pressure terms omitted).

## 6. Search for Bäcklund transformations using Lie formalism.

In some particular cases, it is possible to establish a mapping from solutions of one PDE into the set of solutions of another, i.e. to find a Bäcklund transformation, using Lie formalism.

Suppose we are given a set of differential equations $\mathbf{E}(\mathbf{u}, \mathbf{x}) = 0$, and $\mathbf{u} = \mathbf{F}(\mathbf{x})$ is a solution. Let us, as in (1.1), transform the solution:

$$x_1^i = f^i(\mathbf{x}, \mathbf{u}, \mathbf{u}_1, a),$$   (6.1)

$$u_1^k = g^k(\mathbf{x}, \mathbf{u}, \mathbf{u}_1, a),$$

but instead of $\underset{1}{X}\, \mathbf{E}(\mathbf{u},\, \mathbf{x})\big|_{\mathbf{E}(\mathbf{u},\, \mathbf{x})=0} = 0$ demand

$$\underset{1}{X}\, \mathbf{E}(\mathbf{u},\, \mathbf{x})\big|_{\mathbf{E}(\mathbf{u},\, \mathbf{x})=0} = F(\mathbf{x}, \mathbf{u}, \mathbf{u}_1),$$   (6.2)

where F is some function.

Then, after solving (6.2) and using Lie theorem to reconstruct the transformation, we will get a mapping of solutions of $\mathbf{E}(\mathbf{u}, \mathbf{x}) = 0$ into solutions of a *different equation*.

First, we made an attempt to find a transition from pure magnetic field (4.1) to force-free configurations (5.1), with an assumption $\alpha = \alpha(\mathbf{u}, \mathbf{u}_1)$. The result is that only $\alpha = 0$ is admissible.

Second, we searched for transitions from force-free configurations represented in the form

$$\mathrm{curl}(\mathbf{B}) \times \mathbf{B} = 0,$$   (6.3)
$$\mathrm{div}\, \mathbf{B} = 0$$

to plasma equilibrium equations (2.1), where the unknown pressure $P$ was assumed to depend on $(\mathbf{x}, \mathbf{u}, \mathbf{u}_1)$. The only solution appeared to be $P$ = const, therefore it is impossible to build a Bäcklund transformation from (6.3) to (2.1) using this method.

## *References.*